\input amssym.def
\input epsf

\let \blskip = \baselineskip
\parskip=1.2ex plus .2ex minus .1ex

\tabskip 20pt
\tolerance = 1000
\pretolerance = 50
\newcount\itemnum
\itemnum = 0
\overfullrule = 0pt

\def\title#1{\bigskip\centerline{\bigbigbf#1}}
\def\author#1{\bigskip\centerline{\bf #1}\smallskip}
\def\address#1{\centerline{\it#1}}
\def\abstract#1{\vskip1truecm{\narrower\noindent{\bf Abstract.} #1\bigskip}}

\def\sp{\bigskip}
\def\nosp{\vskip -\the\blskip plus 1pt minus 1pt}

\def\br{\hfil\break} 
\def\ti{\br \hglue \the \parindent}

\def\ce#1{\LP\centerline{#1}}

\def\skipit#1{}
\def\mdag{\raise 3pt\hbox{\dag}}

\def\XP{\par\noindent\hang}
\def\LP{\par\noindent}
\def\BP[#1]{\par\item{[#1]}}
\def\SH#1{\sp\vskip\parskip\leftline{\bigbf #1}\nobreak}

\def\TH#1{\sp\XP{\bf THEOREM\ \shead#1}}

\def\PF{\LP{\bf Proof:\ }}
\def\NX{\advance\itemnum by 1 \sp\LP {\bf \shead \the\itemnum.\ }}
\def\qed{\null\nobreak\hfill\hbox{${\vrule width 5pt height 6pt}$}\par\sp}

\def\cart{\>\hbox{${\vcenter{\vbox{
    \hrule height 0.4pt\hbox{\vrule width 0.4pt height 4.5pt
    \kern4pt\vrule width 0.4pt}\hrule height 0.4pt}}}$}\>}
\def\bxmu{\>\hbox{${\vcenter{\vbox {
    \hrule height 0.4pt\hbox{\vrule width 0.4pt height 4pt
    \hskip -1.3pt\lower 1.8pt\hbox{$\times$}\negthinspace\vrule width 0.4pt}
    \hrule height 0.4pt}}}$}\>}

\def\lin#1{\hbox to #1true in{\hrulefill}}


\def\foot#1{\raise 6pt \hbox{#1} \kern -3pt}

\def\fig #1 #2 #3 #4 #5 {\sp \ce{ {\epsfbox[#1 #2 #3 #4]{figs/#5.ps}} }}



\def\JGT{{\it J.\ Graph Theory}}

\def\COMB{{\it Combinatorica}}
\def\DM{{\it Discrete Math.{}}}

\def\CMSJB{{\it Coll.\ Math.\ Soc.\ J\'anos Bolyai}}

			
 \def\th{\theta}





 \def\Imp{\Rightarrow} 
  

\def\({\left(}	\def\){\right)}


\def\CH#1#2{{{#1}\choose{#2}}}

\def\CL#1{\left\lceil{#1}\right\rceil}


\def\SET#1:#2{\{#1\colon\;#2\}}

		
\def\C#1{\left | #1 \right |}    






\magnification=\magstep1
\vsize=9.0 true in
\hsize=6.5 true in
\headline={\hfil\ifnum\pageno=1\else\folio\fi\hfil}
\footline={\hfil\ifnum\pageno=1\folio\else\fi\hfil}

\parindent=20pt
\baselineskip=12pt
\parskip=.5ex  

\def\shead{ }

\font\bigbf = cmb10 scaled \magstep1

\font\bigbigbf = cmb10 scaled \magstep2


\def\gpic#1{#1
     \bigskip\par\noindent{\centerline{\box\graph}}
     \medskip} 

\def\title#1{\bigskip\centerline{\bigbigbf#1}}
\def\author#1{\bigskip\centerline{\bf #1}}
\def\address#1{\centerline{\it#1}}
\def\abstract#1{\vskip10pt{\narrower\noindent{\bf Abstract.} #1\bigskip}}

\title{STAR-FACTORS OF TOURNAMENTS}
\author{Guantao Chen\foot{*}}
\address{Georgia State University, Atlanta, GA 30303}
\author{Xiaoyun Lu}
\address{Chinese University of Hong Kong, Hong Kong}
\author{Douglas B. West\foot{\dag}}
\address{University of Illinois, Urbana, IL 61801-2975}
\vfootnote{}{\br
   \foot{*}Research supported in part by NSA/MSP Grant MDA904-94-H-2060.\br
   \foot{\dag}Research supported in part by NSA/MSP Grant MDA904-93-H-3040.\br
   Running head: STAR-FACTORS OF TOURNAMENTS\br
   Keywords: Tournament, spanning subgraph, star\br
   Written June 1995; revised June 1997.}
\abstract{
   Let $S_m$ denote the $m$-vertex simple digraph formed by $m-1$
   edges with a common tail.  Let $f(m)$ denote the minimum $n$ such that every
   $n$-vertex tournament has a spanning subgraph consisting of $n/m$
   disjoint copies of $S_m$.  We prove that
   $m\lg m-m\lg \lg m \le f(m)\le 4m^2-6m$ for sufficiently large $m$.}

If $D$ is an acyclic digraph of order $n$, then $D$ occurs as a subgraph of
every tournament with at least $2^{n-1}$ vertices; this follows easily by
induction on $n$.  Special digraphs can be guaranteed to appear even when
there are no extra vertices.  A digraph of order $n$ is {\it unavoidable} if it
appears in every $n$-vertex tournament.  A {\it claw} is a digraph obtained by
identifying the sources of a set of edge-disjoint paths.  The question of which
$n$-vertex claws are unavoidable is studied in [5,6,10].  Further references and
additional unavoidable directed trees appear in [7].

A claw formed using paths of length one is a {\it star}; the $m$-{\it star}
$S_m$ consists of $m-1$ edges with a common tail and $m-1$ distinct heads.  We
use $kS_m$ to denote a disjoint union of $k$ stars of order $m$.  A 3-vertex
tournament need not contain a spanning $S_3$, but every 6-vertex tournament
contains a spanning $2 S_3$.  Some 8-vertex tournaments avoid $2 S_4$ (see
Theorem 3).  A copy of $kS_m$ in a tournament of order $km$ is an
{\it $m$-star-factor} of the tournament.

If $m$-star-factors are unavoidable for some multiple of $m$, then inductively
they are unavoidable for all larger multiples of $m$.  Let $f(m)$ be the least
$n$ such that every $n$-vertex tournament contains an $m$-star-factor (if such
$n$ exists).  We prove that $m\lg m-m\lg \lg m \le f(m)\le 4m^2-6m$ for
sufficiently large $m$.  The upper bound holds for all $m$.  (Always $\lg$
denotes $\log_2$.)

Reid [9] asked for the minimum $n$ such that all $n$-vertex tournaments have $k$
pairwise-disjoint transitive subtournaments of order $m$.  For fixed $m$,
Erd\H os proved the existence of $g(m)$ such that vertices of tournaments with
order $g(m)$ can be partitioned into sets inducing transitive tournaments of
order $m$.  Since every transitive tournament has a spanning star,
$f(m)\le g(m)$, and our lower bound requires $g(m)\ge m\lg m-m\lg\lg m$ (see our
final remark for an upper bound on $g(m)$).  (Lonc and Truszcy\'nski [4]
observed the weaker conclusion that the vertices of a sufficiently large
tournament can be partitioned into sets of size {\it at least} $k$ that induce
transitive subtournaments, but their result holds in a more general setting.) 

Our lower bound depends on a bound for another problem.  We say that a vertex of
a tournament {\it dominates} all its successors, and a tournament {\it is
$k$-dominated} (historically, ``{\it satisfies property $P_k$}'') if for every
$k$-set $U$ of vertices there is a vertex outside $U$ that dominates $U$.
Erd\H os [2] proved that when $n$ is sufficiently large, some $n$-vertex
tournament is $k$-dominated, because then the expected number of undominated
$k$-sets in the random tournament is less than 1.  The expected number of
undominated $k$-sets is $\CH nk(1-2^{-k})^{n-k}$.  Since
$\CH nk(1-x)^t<(ne/k)^ke^{-xt}$, the
expectation is less than 1 when $n>2^k k^2 (\ln 2)(1+o(1))$.
(Erd\H os and Moser [3] obtained even stronger conclusions when $n$ is above
this threshold.)

Let $h(k)$ be the minimum
$n$ such that some $n$-vertex tournament is $k$-dominated.  The bounds are
$ck2^k\le h(k)\le2^k k^2 (\ln 2)(1+o(1))$; we will use the upper bound.
(The lower bound is by Szekeres and described in [8]; Erd\H os had proved that
$h(k)\ge 2^{k+1}-1$.  The upper bound is that of Erd\H os described above; see
[1, p.5-6] for further discussion of this problem and these computations).  
\TH 1.
If $n< m\lg m-m\lg\lg m$ (for sufficiently large $m$), then some $n$-vertex
tournament has no spanning $m$-star-factor.
\PF
We show first that $kS_m$ is avoidable when $km\ge h(k)$.  Since $km\ge h(k)$,
some $km$-vertex tournament $T$ is $k$-dominated.  Let $U$ be a set of $k$
sources of stars in $T$.  Since $T$ is $k$-dominated, some vertex in $T$
dominates $U$.  Such a vertex belongs to none of these $k$ stars.  Thus $T$ is
not spanned by $k$ stars of any sizes.

Thus it suffices to show that $n\ge h(n/m)$ when $n$ is a multiple of $m$ such
that $n<m(\lg m-\lg\lg m)$.  Let $k=n/m$.  Since
$h(k)\le2^k k^2 (\ln 2)(1+o(1))$, it suffices to show that
$n\ge 2^{(n/m)}(n/m)^2(\ln 2)(1+o(1))$.  This simplifies to
$2m\lg m \ge n+m(\lg n+O(1))$, which is satisfied for
$n<m(\lg m - \lg \lg m)$ when $m$ is sufficiently large.  \qed

\vskip -1pc
\TH 2.
If $n>4m^2-6m$ and $n$ is a multiple of $m$, then every $n$-vertex tournament
has an $m$-star-factor.

\PF
Consider an $n$-vertex tournament, and let $x$ be a vertex of maximum
out-degree.  We construct an $m$-star-factor.  In the subtournament induced by
$N^-(x)$, we use as many disjoint $m$-stars as possible.  Let $A$ be the subset
of $N^-(x)$ not covered by these stars.  We have $\C A \le 2m-3$, else within
$A$ we can find another $m$-star.  If some $m$-star has source in $A$ and
leaves in $N^+(x)$, we use it.  When no further $m$-stars of this type can be
found, let $A'$ be the remaining subset of $N^-(x)$, and let $a=\C{A'}$.
Figure 1 illustrates these sets.

Let $B$ be the remaining subset of $N^+(x)$, and let $B'$ be the set of
vertices in $B$ that dominate $A'$.  By the definition of $A'$, there are at
most $a(m-2)$ edges from $A'$ to $B$.  At least $\C B-a(m-2)$ vertices of $B$
are untouched by these edges, and hence $\C{B'}\ge\C B - a(m-2)$.  Since
$d^+(x)\ge (n-1)/2$ and at most $2m-3-a$ stars were used with source in $A$ and
leaves in $N^+(x)$, we have $\C B\ge (n-1)/2-(m-1)(2m-3-a)$.  Hence
$\C{B'}\ge(n-1)/2-(m-1)(2m-3)+a$.

Since $n>2m(2m-3)$, we have $\C{B'}\ge2m-3$.  This means that the subtournament
induced by $B'$ has a vertex $y$ with out-degree at least $m-2$.  Let $R$ be a
subset of $B'$ consisting of $y$ and $m-2$ successors of $y$.  If $a\ge m$,
we use a vertex of $B'$ outside of $R$ as the source of an $m$-star dominating
$m-1$ vertices in $A'$.  Since $a\le 2m-3$, fewer than $m$ vertices of $A'$
remain uncovered.  If any remain uncovered, we use $y$ as the source of an
$m$-star that dominates the rest of $A'$ (or all of $A'$ if $0<a<m$) and as much
of $R$ as needed to complete the star.

We have now covered all vertices except $x$ and a subset of $B$.  We iteratively
use $m$-stars from the subtournament induced by $B$.  As long as at least $2m-2$
vertices remain in $B$, we can find another $m$-star.  When fewer than $2m-2$
vertices remain, we use a star with source $x$.  Now fewer than $m-1$ vertices
remain, which must equal 0 since $n$ is a multiple of $m$.  \qed

\gpic{
\expandafter\ifx\csname graph\endcsname\relax \csname newbox\endcsname\graph\fi
\expandafter\ifx\csname graphtemp\endcsname\relax \csname newdimen\endcsname\graphtemp\fi
\setbox\graph=\vtop{\vskip 0pt\hbox{%
    \special{pn 8}%
    \special{ar 1171 805 366 732 0 6.28319}%
    \special{ar 2634 805 366 732 0 6.28319}%
    \special{pa 912 1322}%
    \special{pa 1429 287}%
    \special{fp}%
    \special{pa 2375 1322}%
    \special{pa 2893 287}%
    \special{fp}%
    \graphtemp=.5ex\advance\graphtemp by 0.000in
    \rlap{\kern 1.171in\lower\graphtemp\hbox to 0pt{\hss $N^-(x)$\hss}}%
    \graphtemp=.5ex\advance\graphtemp by 0.000in
    \rlap{\kern 2.634in\lower\graphtemp\hbox to 0pt{\hss $N^+(x)$\hss}}%
    \graphtemp=.5ex\advance\graphtemp by 0.930in
    \rlap{\kern 1.277in\lower\graphtemp\hbox to 0pt{\hss $A$\hss}}%
    \graphtemp=.5ex\advance\graphtemp by 0.930in
    \rlap{\kern 2.741in\lower\graphtemp\hbox to 0pt{\hss $B$\hss}}%
    \graphtemp=.5ex\advance\graphtemp by 1.427in
    \rlap{\kern 1.171in\lower\graphtemp\hbox to 0pt{\hss $A'$\hss}}%
    \graphtemp=.5ex\advance\graphtemp by 1.427in
    \rlap{\kern 2.634in\lower\graphtemp\hbox to 0pt{\hss $B'$\hss}}%
    \graphtemp=.5ex\advance\graphtemp by 0.805in
    \rlap{\kern 0.073in\lower\graphtemp\hbox to 0pt{\hss $|A|\le2m-3$\hss}}%
    \graphtemp=.5ex\advance\graphtemp by 1.171in
    \rlap{\kern 0.073in\lower\graphtemp\hbox to 0pt{\hss $|A'|=a$\hss}}%
    \special{ar 1171 1723 439 439 -2.143215 -0.998378}%
    \special{ar 2634 1723 439 439 -2.143215 -0.998378}%
    \graphtemp=.5ex\advance\graphtemp by 0.220in
    \rlap{\kern 1.079in\lower\graphtemp\hbox to 0pt{\hss $\bullet$\hss}}%
    \graphtemp=.5ex\advance\graphtemp by 0.530in
    \rlap{\kern 0.933in\lower\graphtemp\hbox to 0pt{\hss $\bullet$\hss}}%
    \graphtemp=.5ex\advance\graphtemp by 0.530in
    \rlap{\kern 1.079in\lower\graphtemp\hbox to 0pt{\hss $\bullet$\hss}}%
    \graphtemp=.5ex\advance\graphtemp by 0.530in
    \rlap{\kern 1.226in\lower\graphtemp\hbox to 0pt{\hss $\bullet$\hss}}%
    \special{pa 1079 220}%
    \special{pa 1167 406}%
    \special{fp}%
    \special{sh 1.000}%
    \special{pa 1152 332}%
    \special{pa 1167 406}%
    \special{pa 1119 348}%
    \special{pa 1152 332}%
    \special{fp}%
    \special{pa 1167 406}%
    \special{pa 1226 530}%
    \special{fp}%
    \special{pa 1079 220}%
    \special{pa 991 406}%
    \special{fp}%
    \special{sh 1.000}%
    \special{pa 1039 348}%
    \special{pa 991 406}%
    \special{pa 1006 332}%
    \special{pa 1039 348}%
    \special{fp}%
    \special{pa 991 406}%
    \special{pa 933 530}%
    \special{fp}%
    \special{pa 1079 220}%
    \special{pa 1079 406}%
    \special{fp}%
    \special{sh 1.000}%
    \special{pa 1098 333}%
    \special{pa 1079 406}%
    \special{pa 1061 333}%
    \special{pa 1098 333}%
    \special{fp}%
    \special{pa 1079 406}%
    \special{pa 1079 530}%
    \special{fp}%
    \graphtemp=.5ex\advance\graphtemp by 0.732in
    \rlap{\kern 1.445in\lower\graphtemp\hbox to 0pt{\hss $\bullet$\hss}}%
    \graphtemp=.5ex\advance\graphtemp by 0.878in
    \rlap{\kern 2.360in\lower\graphtemp\hbox to 0pt{\hss $\bullet$\hss}}%
    \graphtemp=.5ex\advance\graphtemp by 0.732in
    \rlap{\kern 2.360in\lower\graphtemp\hbox to 0pt{\hss $\bullet$\hss}}%
    \graphtemp=.5ex\advance\graphtemp by 0.585in
    \rlap{\kern 2.360in\lower\graphtemp\hbox to 0pt{\hss $\bullet$\hss}}%
    \special{pa 1445 732}%
    \special{pa 1994 644}%
    \special{fp}%
    \special{sh 1.000}%
    \special{pa 1919 637}%
    \special{pa 1994 644}%
    \special{pa 1925 674}%
    \special{pa 1919 637}%
    \special{fp}%
    \special{pa 1994 644}%
    \special{pa 2360 585}%
    \special{fp}%
    \special{pa 1445 732}%
    \special{pa 1994 820}%
    \special{fp}%
    \special{sh 1.000}%
    \special{pa 1925 790}%
    \special{pa 1994 820}%
    \special{pa 1919 826}%
    \special{pa 1925 790}%
    \special{fp}%
    \special{pa 1994 820}%
    \special{pa 2360 878}%
    \special{fp}%
    \special{pa 1445 732}%
    \special{pa 1994 732}%
    \special{fp}%
    \special{sh 1.000}%
    \special{pa 1921 713}%
    \special{pa 1994 732}%
    \special{pa 1921 750}%
    \special{pa 1921 713}%
    \special{fp}%
    \special{pa 1994 732}%
    \special{pa 2360 732}%
    \special{fp}%
    \graphtemp=.5ex\advance\graphtemp by 0.073in
    \rlap{\kern 1.902in\lower\graphtemp\hbox to 0pt{\hss $\bullet$\hss}}%
    \graphtemp=.5ex\advance\graphtemp by 0.000in
    \rlap{\kern 1.902in\lower\graphtemp\hbox to 0pt{\hss $x$\hss}}%
    \special{pa 1482 384}%
    \special{pa 1734 198}%
    \special{fp}%
    \special{sh 1.000}%
    \special{pa 1560 270}%
    \special{pa 1734 198}%
    \special{pa 1614 343}%
    \special{pa 1560 270}%
    \special{fp}%
    \special{pa 1734 198}%
    \special{pa 1902 73}%
    \special{fp}%
    \special{pa 1902 73}%
    \special{pa 2155 260}%
    \special{fp}%
    \special{sh 1.000}%
    \special{pa 2035 114}%
    \special{pa 2155 260}%
    \special{pa 1981 188}%
    \special{pa 2035 114}%
    \special{fp}%
    \special{pa 2155 260}%
    \special{pa 2323 384}%
    \special{fp}%
    \special{pa 2451 1445}%
    \special{pa 1793 1445}%
    \special{fp}%
    \special{sh 1.000}%
    \special{pa 1976 1491}%
    \special{pa 1793 1445}%
    \special{pa 1976 1399}%
    \special{pa 1976 1491}%
    \special{fp}%
    \special{pa 1793 1445}%
    \special{pa 1354 1445}%
    \special{fp}%
    \special{pn 28}%
    \special{pa 1482 384}%
    \special{pa 1902 73}%
    \special{fp}%
    \special{pa 1902 73}%
    \special{pa 2323 384}%
    \special{fp}%
    \special{pa 2451 1445}%
    \special{pa 1354 1445}%
    \special{fp}%
    \hbox{\vrule depth1.573in width0pt height 0pt}%
    \kern 3.000in
  }%
}%
}
\sp
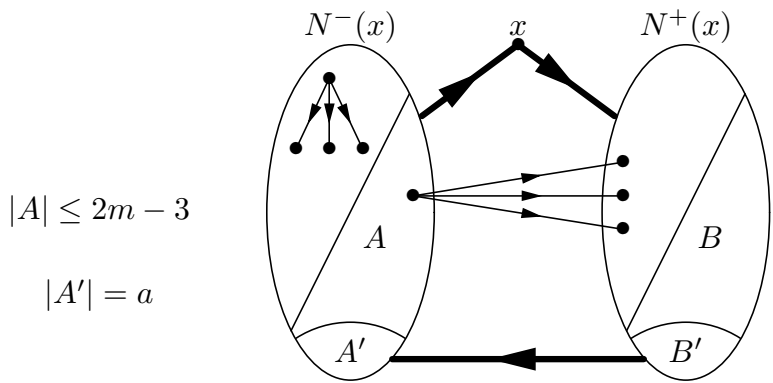
\ce{Figure 1.  Structure of the $m$-star factor in Theorem 2.}
\sp

For the interested reader, we discuss the first few values of $f(m)$.

\TH{3.}
$f(2)=2$, $f(3)=6$, and $f(4)\ge12$.
\PF
Trivially, $f(2)=2$.  The cyclic triple has no $S_3$, so $f(3)\ge6$.
Let $T$ be an arbitrary tournament of order 6.  Let $a|bc$ denote a
star centered at $a$ with the remaining vertices as leaves, and let $(abc)$
denote a cycle with edges $ab,bc,ca$.  Every subtournament induced by at least
four vertices contains $S_3$, so we begin with $x|yz$ and the edge $yz$.  Given
any $S_3$ in $T$, we are finished unless the remaining three vertices form a
cyclic triple, so we also have $(uvw)$.  If $d^+(x)=5$, then we have a 3-star
within $T-x$ and another with center $x$.  Thus we may assume the edge $ux$.
Now $u|xv \Imp (yzw)$, $w|uy \Imp (xzv)$, and $v|xw \Imp (yzu)$.
We now have $2S_3$ formed by $z|vw$ and $u|xy$.

For $m=4$, we present 8-vertex tournaments having no $2S_4$.  Our examples
contain a particular 6-vertex tournament $T_6$.  Construct the 3-cycles
$(x_1x_2x_3)$ and $(y_3y_2y_1)$.  Add the three edges of the form $x_iy_i$ and
all six edges of the form $y_ix_j$ for $i\ne j$.  We have $d^+(x_i)=2$ and
$d^+(y_i)=3$.  Each copy of $S_4$ in this tournament has $y_i$ as its center,
for some $i$, with leaves $y_{i-1},x_{i-1},x_{i+1}$ (indices modulo 3).  Let
$X=\{x_1,x_2,x_3\}$ and $Y=\{y_1,y_2,y_3\}$.

Add a vertex $u$ such that $N^-(u)=X$ and $N^+(u)=Y$.  The resulting 7-vertex
tournament $T_7$ is regular and does not contain $S_4+S_3$.  Each copy of
$S_4$ consists of a vertex and all its successors.  Deleting the 4-star
centered at $u,x_i,y_i$ leaves the 3-cycles $(x_1x_2x_3)$,
$(x_{i-1}y_{i-1}y_{i+1})$, $(y_{i+1}x_iu)$, respectively.

To construct an 8-vertex tournament avoiding $2S_4$, it suffices to add a sink
to $T_7$.  An 8-vertex tournament with a sink contains $2S_4$ if and only if the
7-vertex tournament obtained by deleting the sink contains $S_4+S_3$, which
$T_7$ does not.  Another such tournament, with no sink, is obtained by adding
to $T_7$ a vertex $v$ such that $N^-(v)=X$ and $N^+(v)=Y\cup\{u\}$.  In this
tournament $T_8$, the vertices of $Y\cup\{u\}$ have out-degree 3, and those
of $X\cup\{v\}$ have out-degree 4.  Deleting the 4-star centered at $u$ leaves
a 4-vertex tournament with a sink.  Deleting the 4-star centered at $y_i$
leaves $x_i,v$ with out-degree 2 and $y_{i+1},u$ with out-degree 1.
Thus the centers of a copy of $2S_4$ in $T_8$ must both lie in $X\cup\{v\}$.
Two centers from $X$ cannot dominate all of $Y$, and a center from $X$ along
with $v$ as a second center cannot dominate all of $X$.  \qed

One referee observed that $f(3)=6$ also follows immediately from the first
sentence of [9].  For $f(4)$, Theorem 2 yields $f(4)\le 24$.  We convinced
ourselves that $f(4)=12$ via a case analysis too tedious and painful to
reproduce.  We leave this open in the hope that someone will improve the
general bound in Theorem 2.


Finally, we thank Zbigniew Lonc (Warsaw University of Technology) for providing
a proof of a bound on $g(m)$, the existence of which was attributed to Erd\H os
in [9] but appears never to have been published.  First, let $\th(m)$ denote the
minimum number of vertices that forces every tournament to have a transitive
subtournament of order $m$.  The bound $\th(m)\le2^{m-1}$ is usually attributed
to Erd\H os and Moser and appears in Moon [8, p.15].

Lonc argued as follows that $g(m)\le m\CL{\th((m-1)\th(m))/m}$.  Let $T$
be a tournament of this order; $T$ has a transitive subtournament $A$ of order
$(m-1)\th(m)$.  From $T-A$ we delete transitive tournaments of order $m$ until
fewer than $\th(m)$ vertices remain outside $A$.  For each remaining vertex
$v\in T-A$, iteratively, we find $m-1$ successors or $m-1$ predecessors of $v$
in $A$.  Together with $v$, these form a transitive subtournament.  When we
process the last such vertex outside $A$, there remain at least $2(m-1)$
vertices in $A$, so we can complete this phase.  Finally, since $A$ is
transitive, the remaining vertices of $A$ are partitioned into transitive
subtournaments of order $m$.

\SH
{\ce{References}}
\frenchspacing
\BP[1]
N.~Alon and J.H.~Spencer, {\it The Probabilistic Method}, (Wiley 1992).
\BP[2]
P. Erd\H os, On a problem in graph theory, {\it Math. Gaz.} 47(1963), 220--223.
\BP[3]
P. Erd\H os and L. Moser, A problem on tournaments,
{\it Canad. Math. Bull.} 7(1964), 351--356.
\BP[4]
Z. Lonc and M. Truszcy\'nski, Decomposition of large uniform hypergraphs,
{\it Order} 1(1985), 345--350.
\BP[5]
X.~Lu, On claws belonging to every tournament, \COMB\ 11(1991), 173--179.
\BP[6]
X.~Lu, Claws contained in all $n$-tournaments, \DM\ 119(1993), 107-111.
\BP[7]
X.~Lu, On avoidable and unavoidable trees, \JGT\ 22(1996), 335--346.
\BP[8]
J.W.~Moon, {\it Topics on Tournaments}, (Holt, Reinhart \& Winston 1968).
\BP[9]
K.B.~Reid, Three problems on tournaments, {\it Graph Theory and its
Applications: East and West (Proc.\ 2nd China-USA Intl.\ Conf.\ Graph Th.,
San Francisco 1989)}, New York Acad.\ Sci. {\it Annals} 576(1989), 466--473.
\BP[10]
M.~Saks and V.T.~S\'os, On unavoidable subgraphs of tournaments,
in {\it Finite and Infinite Sets (Eger),} \CMSJB\ 37(1981), 663--674.
\bye